\documentclass[11pt]{article}
\usepackage{color}
\usepackage{amsfonts}	
\newtheorem{theo}{Theorem}[section]
\newtheorem{lem}[theo]{Lemma}
\newtheorem{prop}[theo]{Proposition}

\newcommand{\mysection}[1]{\section{#1} \setcounter{equation}{0}}
\newcommand{\dfrac}[2]{\displaystyle \frac{#1}{#2}}

\newcommand{\proof}{{\sc Proof.} \quad}

\newcommand{\R}{\mathbb{R}}

\newcommand{\be}{\begin{equation} \label}
\newcommand{\ee}{\end{equation}}
\newcommand{\bes}{\begin{equation} \begin{array}{c} \label}
\newcommand{\ees}{\end{array} \end{equation}}
\newcommand{\bea}{\begin{eqnarray}\label}
\newcommand{\eea}{\end{eqnarray}}
\newcommand{\beas}{\begin{eqnarray} \begin{array}{rcl} \label}
\newcommand{\eeas}{\end{array} \end{eqnarray}}
\newcommand{\bas}{\begin{eqnarray*}}\newcommand{\eas}{\end{eqnarray*}}
\newcommand{\bass}{\begin{eqnarray*} \begin{array}{rcl}}
\newcommand{\eass}{\end{array} \end{eqnarray*}}
\newcommand{\basss}{\begin{eqnarray*} \begin{array}{c}}
\newcommand{\easss}{\end{array} \end{eqnarray*}}

\newcommand{\bit}{\begin{itemize}}
\newcommand{\eit}{\end{itemize}}
\newcommand{\nn}{\nonumber}
\newcommand{\eps}{\varepsilon}
\newcommand{\abs}{\\[3mm]}
\newcommand{\parab}{{\cal{P}}}
\newcommand{\opa}{{\cal{A}}}
\newcommand{\upsi}{\underline{\psi}}
\newcommand{\opsi}{\overline{\psi}}
\newcommand{\opsio}{\overline{\psi}_{out}}
\newcommand{\opsii}{\overline{\psi}_{in}}
\newcommand{\uv}{\underline{v}}
\newcommand{\ov}{\overline{v}}
\newcommand{\ren}{\mathbb{R}^n}
\def\qed{\,\unskip\kern 6pt \penalty 500
\raise -2pt\hbox{\vrule \vbox to8pt{\hrule width 6pt
\vfill\hrule}\vrule}\par}
\begin{document}
\title{\bf A Continuum of Extinction Rates\\
for the Fast Diffusion Equation}
\author{\Large
Marek Fila\footnote{fila@fmph.uniba.sk}  \\[4pt]
\Large Juan Luis V\'azquez\footnote{juanluis.vazquez@uam.es} \\[4pt]
\Large Michael Winkler\footnote{michael.winkler@uni-due.de}
}
\date{}
\maketitle

\begin{abstract}
\noindent We find a continuum of extinction rates for solutions $u(y,\tau)\ge 0$ of the fast diffusion equation $u_\tau=\Delta u^m$ in a subrange of exponents $m\in (0,1)$. The equation is posed in $\ren$ for times up to the extinction time $T>0$. The rates take the form $\|u(\cdot,\tau)\|_\infty\sim (T-\tau)^\theta$ \ for a whole interval of  $\theta>0$. These extinction rates depend explicitly on the spatial decay rates of
initial data.

\end{abstract}
\normalcolor
\mysection{Introduction}

We consider the Cauchy problem for the fast diffusion equation:
\begin{equation}\label{fd}
        \left\{ \begin{array}{ll}
        u_\tau = \Delta (u^m/m),
                \qquad & y\in\R^n, \ \tau\in (0,T), \\[2mm]
        u(y,0)=u_0(y)\ge 0, \qquad & y\in\R^n,
        \end{array} \right.
\end{equation}
where $m\in (0,1)$ and $T>0$.  The factor $1/m$ is not essential; it is inserted
into the equation for normalization so that it  can also be written as
$ u_ \tau=\nabla\cdot(u^{m-1}\,\nabla u)$. In that way, it is readily seen that
the diffusion coefficient $c(u)=u^{m-1}\to\infty$ as $u\to 0$ if $m<1$, hence the
name Fast Diffusion Equation  (but notice that $c(u)\to 0$ as $u\to\infty$).
Furthermore, it is known that for $m$ below a critical exponent $m_c=(n-2)/n$ all
solutions with initial data in some convenient space, like $L^p(\ren)$ with $p=n(1-m)/2$,
extinguish in finite time. We will always work in this range, $m<m_c$, and consider
solutions which vanish in a finite time.  The purpose of this paper is to study the
rates of extinction of such solutions. Our main contribution  is to provide a continuum
of rates of extinction for fixed $m$. Technical reasons  imply that $m$  must be in the
range $0<m<m_*=(n-4)/(n-2)$, $n\ge 5$, for the construction to work. This restriction may be essential.

 Let us review the state of the question from a broader perspective.
The description of the asymptotic behaviour of  the global in time  solutions
of (\ref{fd}) as $\tau\to\infty$ for $m\ge m_c$ is a very active subject,
and the study has been extended in recent times to the behaviour near extinction
for $m<m_c$, both in bounded domains or in the whole space.  In the former case,
the rate of decay for bounded solutions  is universal, of the form
$\|u(\cdot,\tau)\|_\infty = O((T-\tau)^{1/(1-m)})$ when $m> m_s=(n-2)/(n+2)$,
cf. \cite{BH, FS2000},  but the question is more complicated when $m\le m_s$.

 In the case of the whole space, which is the one of interest  here, the book \cite{Vsmooth} contains a general description of the phenomenon of extinction, where it is explained that not only the occurrence of extinction depends on the size of the initial data, but also that  different initial data may give rise to  different extinction rates, even for the same extinction time; this may happen for all $0<m<m_c$.  It is also proved, cf. \cite{Vsmooth} and quoted references, that the size of the initial data at infinity (the tail of $u_0$) is very important in determining both the extinction time and the decay rates.

Special attention has been given recently to particular classes of data that produce definite estimates. This happens in the case of data with the maximal decay rate compatible with extinction in finite time, which is
\begin{equation}
 u_0(y)\sim A|y|^{-\mu}, \quad \mu:=2/(1-m)
 \end{equation}
as $|y|\to\infty$. Note that $\mu<n$ for $m<m_c$ so these data are not integrable. \normalcolor Thus, the papers \cite{McD, BBDGV, PNAS}
 are concerned with the stabilization as $\tau\to T$ of general solutions towards some special self-similar solutions $U_{D, T}$ known as the 
\emph{generalized Barenblatt solutions}, given by the formula
\begin{equation}\label{baren.form1}
U_{D, T}(y,\tau):=\frac 1{R(\tau)^n} \left(D+\frac{\beta(1-m)}{2}
\left|\frac{y}{R(\tau)}\right|^2\right)^{-\frac{1}{1-m}},
\end{equation}
where for $m<m_c$ we put $R(\tau):=(T-\tau)^{-\beta},$ and
\[
\beta:= \frac1{n(1-m)-2}=\frac{1}{n\,(m_c-m)}=\frac{\mu}{2(n-\mu)}\,.
\]
Here $T\ge 0$ (extinction time) and $D>0$ are free parameters. Note that $R$ depends on $T$.
 It has been proved that the corresponding Barenblatt solutions with exponent $m>m_c$ play 
the role of the Gaussian solution of the linear diffusion equation in describing 
the asymptotic behaviour of a very wide class of nonnegative solutions, i.e., those 
with initial data in $L^1(\ren)$, cf. \cite{Vbook}. To some extent, the solutions 
(\ref{baren.form1}) play a similar role for $m<m_c$ but their basin of attraction 
may be much smaller. This is precisely described in \cite{PNAS}, with  results on the basin of attraction of the family of generalized Barenblatt solutions;  it establishes the optimal rates of convergence of the solutions of (\ref{fd}) towards a unique attracting limit state in that family. All of these solutions will have a decay rate near extinction of the form $\|u(\cdot,\tau)\|_\infty = O((T-\tau)^{n\beta})$, and it is clear that $n\beta>1/(1-m)$.

A very interesting limit case occurs if we take $D=0$ in formula (\ref{baren.form1}), and we find the singular solution
\begin{equation}
U_{0, T}(y,\tau):=k_*\,(T-\tau)^{\mu/2}|y|^{-\mu},
\end{equation}
whose attracting properties have not been studied. Note the value $k_*=(2(n-\mu))^{\mu/2}$.

The question that we address here is the following: Can we obtain different decay rates near extinction for bounded data $u_0(y)$ that behave  at infinity in first approximation like the singular solution, i.\,e., $u_0(y)\sim A\,|y|^{-\mu}$ ? We will show that the answer is yes, and actually we will obtain a whole continuum of rates.

\begin{theo}\label{thm-1} Let $u\ge 0$ be a solution of Problem  {\rm (\ref{fd})}, assume that
 \begin{equation}\label{mn}
	n \ge 5 \qquad \mbox{and} \qquad m \in \Big(0\, , \, \frac{n-4}{n-2}\Big),
\end{equation}
and let the initial function $u_0$ be continuous, bounded, and satisfy the conditions:
$$
0\le u_0(y)\le A\,|y|^{-\mu} \quad \mbox{ for all $y\ne 0$ }
$$
and
$$
A\,|y|^{-\mu}-c_1 |y|^{-l}\le u_0(y) \le A\,|y|^{-\mu}-c_2 |y|^{-l} \qquad \mbox{for \ } |y|\ge 1
$$
for some $A,c_1,c_2>0$, and
\begin{equation}\label{eq:l}
\mu+2< l\le L=\mu+\sqrt{2(n-\mu)}.
\end{equation}
Then the solution has complete extinction precisely at the time $T=(A/k_*)^{1-m}>0$, and there are positive constants $K_1, K_2$ such that for $0<\tau<T$ we have
\begin{equation}
K_1(T-\tau)^{\theta} \le \|u(\cdot,\tau)\|_\infty \le K_2(T-\tau)^{\theta},
\end{equation}
where $\theta= \frac{n\mu-\gamma}{2(n-\mu)}>0$, $\gamma=\frac{\mu(l-\mu-2)(n-l)}{l-\mu}$.
\end{theo}

It is easy to check that under the above assumptions $\theta$ covers an interval $[\theta_{min}, \theta_{max})$ with $0<\theta_{min}<\theta_{max}=\mu n/2(n-\mu)=n\beta. $
This is the precise range of extinction rates of these solutions, to be compared with the standard extinction rate $ (T-\tau)^{n\beta}$ of the Barenblatt examples.

As a precedent to this result, the existence of different rates was established in
Theorem~7.4 of \cite{Vsmooth} for all $m<m_c$ by means of the construction of
self-similar solutions of the form $u(y,\tau)=(T-\tau)^{\alpha}f(y\,(T-\tau)^{\beta})$.
In this way a whole interval $(\overline{\alpha},\infty)$ is covered, which extends
the scope of our present theorem. However, $\overline{\alpha}$ (the anomalous exponent)
is not explicit, we obtain only one solution for each time-decay rate and the
dependence of $\alpha$ on the spatial behavior of the data is not analyzed.
Theorem~\ref{thm-1} clarifies these aspects, explaining the delicate relationship
between both limits, $|y|\to \infty$ for $u_0$ and $\tau\to T$ for $u(y,\tau)$.

The proof of the theorem needs techniques that are only natural  after rescaling the problem.  In fact, the rescaled problem allows us to formulate and prove  a more precise result about the dependence of the rate on the tail of the data and the convergence of the spatial shapes. We devote the next section to the presentation of the rescaling transformation, the resulting rescaled equation and the asymptotic convergence plus grow-up result in that context.
Sections~3--5 will be concerned with proving the result for the rescaled
problem. The last section is devoted to comments and open problems.

\medskip

\noindent {\sc Notations. } Throughout the rest of the paper and unless
mention to the contrary, we keep the conditions $n \ge 5 $ and $m < m_*$.
The exponent $m_*$ also plays a big role in the asymptotic results of
\cite{BBDGV, BGV, PNAS}. We also keep the above symbols and variables. In particular, $\mu=2/(1-m)$ so that  $m<m_c$ means $\mu<n$ and $m<m_*$ means $\mu+2<n$.


\mysection{The rescaled flow}

As we have just said, it is very convenient to rescale the flow and rewrite (\ref{fd}) in self-similar variables by introducing  the time-dependent change of variables
\begin{equation}\label{eq:chgvariable}
t:=\frac{1-m}{2}\log\left(\frac{R(\tau)}{R(0)}\right)\quad\mbox{and}\quad
x:=\sqrt{\frac{\beta(1-m)}{2}}\,\frac y{R(\tau)}\,,
\end{equation}
with $R$ as above, and the rescaled function
\begin{equation}\label{eq:chgvariable2}
v(x,t):= R(\tau)^{n}\,u(y,\tau).
\end{equation}
In these new variables, the generalized Barenblatt functions $U_{D, T}(y,\tau)$ are transformed into \emph{generalized Barenblatt profiles} $V_D(x)$, which are stationary:
\begin{equation}\label{newBaren}
V_D(x):=(D+|x|^2)^\frac 1{m-1},\quad x\in \R^n\,.
\end{equation}
If $u$ is a solution to~(\ref{fd}), then $v$
solves the {\sl rescaled fast diffusion equation}
\begin{equation}\label{FPeqn}
v_t=\Delta (v^m/m)+\mu \,\nabla\cdot(x\,v), \quad t> 0\,,\quad x\in \R^n\, ,
\end{equation}
which is a nonlinear Fokker-Planck equation (NLFP). We put
as initial condition $v_0(x):=R(0)^{-n}\,u_0(y)$, where $x$ and~$y$ are related
according to (\ref{eq:chgvariable}) with $\tau=0$, $x=cy$. Roughly speaking,
$v_0$ is a rescaling of $u_0$ depending only on $T$. We have taken the precise
form of this transformation from \cite{PNAS}. Note also that the factors $1/m$ and
$\mu$ in equation (\ref{FPeqn}) can be eliminated by manipulating the change of variables,
but then the expression of the Barenblatt solutions would contain new constants. Thus, in our scaling the singular solution becomes
\begin{equation}
V_0(x)=|x|^{-\mu},\quad x\in \R^n\setminus\{0\}\,.
\end{equation}


\subsection{Main result for the  NLFP equation}

In the following sections we consider the $v$-equation  (\ref{FPeqn}) with initial
data given by a bounded function $0\le v_0\le V_0$, and such that the difference $V_0-v_0$
has a tail controlled by a power rate.  This is our detailed  result about asymptotic
behaviour of the solution whose initial data $v_0(x)$ are perturbations of the steady state $V_0(x)$.

\begin{theo}\label{thm1}
  Assume that $n$ and $m$ are as in {\rm (\ref{mn})}. Suppose that $v_0 $ is continuous, bounded and nonnegative, and fulfils
  \begin{equation}\label{a:v0}
	|x|^{-\mu}-c_1 |x|^{-l}\le v_0(x) \le |x|^{-\mu}-c_2 |x|^{-l} \qquad \mbox{for  } |x|\ge 1,
  \end{equation}
where $l$ is as in {\rm (\ref{eq:l})}
and $c_1,c_2>0$. Assume also that $v_0(x)\le |x|^{-\mu}$ for all $x\ne 0$. Let  $v$ denote the solution  of \ {\rm (\ref{FPeqn})}. Then:

  (i) \ There exist $K_1, K_2>0$ such that for $ t\ge 1$ we have
  \be{convest1}
K_1 \, e^{\gamma \, t} \le	\|v(\cdot,t)\|_\infty \le K_2 \, e^{\gamma\, t},
\qquad \gamma=\gamma(l)=\frac{\mu(l-\mu-2)(n-l)}{l-\mu}.
  \ee
  (ii) \ For each $r_0>$ one can find $C_1, C_2>0$ such that for $t\ge 1$ and $|x|\ge r_0$
  the following holds
  \be{convest2}
C_1 \, e^{-\lambda\, t}\le |x|^{-\mu}-v(x,t) \le  C_2 \, e^{-\lambda\,t}, \qquad
\lambda=\lambda(l)= (l-\mu-2)(n-l).
  \ee
\end{theo}

Let us comment on the contents and scope of the result.

\noindent {\bf 1.} First of all, it states the two main aspects of the convergence
of the solution $v(\cdot,t)$ towards the singular steady state $V_0$:  \ (\ref{convest2})
establishes the uniform convergence of $v(\cdot,t)$ towards $V_0$ in the complement of a ball centered at the origin, with a precise rate that depends explicitly on the tail decay exponent $l$. On the other hand, estimate  (\ref{convest1}) gives the exact rate of growth of the solutions as $t\to \infty$ to account for the approach to the singular value $V_0(0)=+\infty$.

\noindent {\bf 2.} An important feature of the result is the existence of a {\sl continuum of grow-up rates} for $\|v(\cdot,t)\|_\infty$, and a corresponding {\sl continuum of stabilization rates} of $v(\cdot,t)$ towards $V_0$ in the outer region. Note furthermore that as $l$ approaches the lower value $\mu+2$, the rates go to zero. This limit case is on the other hand easier and does not produce any convergence, since we can consider the example of the generalized Barenblatt solutions $V_D$ given in \ (\ref{newBaren}). Indeed, they satisfy $0< V_D< V_0$ and
$$
V_0(x) - V_D(x)= C|x|^{-(\mu+2)}+ o\left(|x|^{-(\mu+2)}\right) \qquad \mbox{as } \ |x|\to \infty.
$$
Since they are stationary, no convergence to $V_0$ holds in this case.

\noindent {\bf 3.} The conditions on $l$ imply that the perturbation $V_0-v_0$ is never integrable,
contrary to the usual assumptions made in variational methods. Let us now examine the maximal
grow-up rate that we have achieved. Note first that
$\gamma(\mu+2)=\gamma(n)=0$. The maximum of $\gamma$ in (\ref{convest1}) is attained at
$l=L$, and
$$
\gamma(L)= \mu\left(n+2-\mu- 2\sqrt{2(n-\mu)}\right).
$$
This is lower than the maximal growth rate of any bounded solution that is given by the growth of the
spatially homogeneous solution $\tilde v(t)=ce^{\mu n t}$.  We conjecture that $\gamma(L)$ is the
largest exponent that can be achieved by the solutions under the conditions of
the theorem, even if we allow $l$ to be larger than $L$. The bound from below follows
immediately from the lower bound in (\ref{convest1}), but to obtain the corresponding
bound from above is still an open problem.

\noindent {\bf 4.} As $m\to m_*$ we have $\mu+2\to n$ and the interval $(\mu+2,L] $
shrinks to the empty set while the admissible values of the exponents $\gamma$ and $\lambda$ go to zero.

\noindent {\bf 5.} Results similar as in Theorem~\ref{thm1}
were obtained for the standard Fujita equation
$$u_t=\Delta u + u^p,\quad x\in \R^n,\quad n>10,\quad
p>\dfrac{(N-2)^2-4N+8\sqrt{N-1}}{(N-2)(N-10)},
$$
in \cite{FKWY, FW, FWY1}.

\noindent {\bf 6.} Finally,  we apply  the results of Theorem~\ref{thm1}~(i)
to prove Theorem~\ref{thm-1}. Notice that under the assumptions of
Theorem~\ref{thm-1}, if we take the prescribed value of $T$ then $v_0$ satisfies the
hypotheses of Theorem~\ref{thm1}, so that the solution $v$ is global in time and
stabilizes to $V_0$; this means that the extinction time of $u$ is precisely $T$.
The extinction rate of $u$ is obtained by rewriting
the bounds in  (\ref{convest1}). Recall that
$$
\|u(\cdot,\tau)\|_\infty=R(\tau)^{-n}\|v(\cdot,t)\|_\infty \sim (T-\tau)^{n\beta}e^{\gamma\, t},
$$
and $T-\tau=T e^{-2(n-\mu)t}$. The conclusion follows.

\mysection{Auxiliary results for the rescaled problem}
After the previous transformation, in the radially symmetric case we end up with the problem
\be{0}
	\left\{ \begin{array}{l}
	\parab v := v_t - \frac{1}{m} \Big( (v^m)_{rr}+\frac{n-1}{r} (v^m)_r \Big) - \mu r v_r - \mu n v \, =0,
		\qquad r>0, \ t>0, \\[2mm]
	v(r,0)=v_0(r), \qquad r \ge 0.
	\end{array} \right.
\ee

An important role is played by the quadratic equation
\be{eqn_alpha}
	\alpha^2-(n-2-\mu-\kappa)\alpha+2\kappa=0,
\ee
where
\be{kappa}
	\kappa=\frac{(l-\mu-2)(n-l)}{l-\mu}
\ee
is positive if $\mu+2<l<n$. The roots $\alpha_-$ and $\alpha_+$ of (\ref{eqn_alpha}) are given by
\be{alpha}
	\alpha_\pm = \frac{n-2-\mu-\kappa \pm \sqrt{(n-2-\mu-\kappa)^2-8\kappa}}{2},
\ee
and the following way to rewrite $\alpha_\pm$ indicates why the value $l=\mu+\sqrt{2(n-\mu)}$ plays an important role
in the sequel (cf.~Section \ref{sect_upper}).
\begin{lem}\label{lem_alpha}
  Assume {\rm (\ref{mn})}, and that $l \in (\mu+2,n)$. Then the roots $\alpha_\pm$
  of {\rm (\ref{eqn_alpha})} can be expressed as follows: \\
  i) \ If $\mu+2<l\le \mu+\sqrt{2(n-\mu)}$ then
  \bas
	\alpha_-=l-\mu-2 \qquad \mbox{and} \qquad \alpha_+=\frac{2(n-l)}{l-\mu}.
  \eas
  ii) \ If $\mu+\sqrt{2(n-\mu)} \le l<n$ then
  \bas
	\alpha_-=\frac{2(n-l)}{l-\mu} \qquad \mbox{and} \qquad \alpha_+=l-\mu-2.
  \eas
\end{lem}
\proof
  Since it can easily be checked that $\alpha_1:=l-\mu-2 < \frac{2(n-l)}{l-\mu}=:\alpha_2$
  if and only if $l < \mu+ \sqrt{2(n-\mu)}$, we only
  need to check that both $\alpha_1$ and $\alpha_2$ solve (\ref{eqn_alpha}). As to $\alpha_1$, this follows from
  \bas
	\alpha_1^2-(n-2-\mu-\kappa) +2\kappa &=& (l-\mu-2)^2-(n-2-\mu-\kappa)(l-\mu-2)+2\kappa \\
	&=& (l-\mu-2)\Big(l-\mu-2-(n-2-\mu-\kappa)\Big)+2\kappa \\
	&=& (l-\mu-2)(l-n) + (l-\mu)\kappa
  \eas
  and the fact that $(l-\mu)\kappa=(l-\mu-2)(n-l)$. Using $\alpha_2+2=\frac{2(n-\mu)}{l-\mu}$, we moreover compute
  \bas
	\alpha_2^2-(n-2-\mu-\kappa)\alpha_2+2\kappa &=& \Big( \frac{2(n-l)}{l-\mu}\Big)^2
	- (n-2-\mu) \frac{2(n-l)}{l-\mu} + (\alpha+2)\kappa \\
	& & \hspace*{-20mm} = \, \frac{2(n-l)}{(l-\mu)^2} \Big( 2(n-l) - (n-2-\mu)(l-\mu)+(n-\mu)(l-\mu-2)\Big),
  \eas
  from which we immediately find that also $\alpha_2$ solves (\ref{eqn_alpha}).
\qed
The following two lemmata apply to parameters $n,m$ and $\kappa$ more general than required in (\ref{mn}) and (\ref{kappa}).
\begin{lem}\label{lem7}
  Let $n\ge 1$, $m>0, \kappa>0$ and $\sigma_0>0$, and set
  \bas
	\sigma(t):= \sigma_0 \, e^{\mu \kappa t}, \qquad t\ge 0,
  \eas
  and
  \bas
	\xi(r,t):=\sigma^\frac{1}{\mu}(t) r, \qquad r\ge 0, \ t\ge 0.
  \eas
  Suppose that $\psi:[0,\infty) \to [0,\infty)$ is twice continuously differentiable in $(\xi_0,\xi_1)$ with some
  $\xi_0$ and $\xi_1$ satisfying $0 \le \xi_0 < \xi_1$. Then for
  \bas
	v(r,t):=\sigma(t) \Big(\xi^2(r,t)+\psi(\xi(r,t)) \Big)^{-\frac{\mu}{2}}, \qquad r\ge 0, \ t\ge 0,
  \eas
  we have the identity
  \be{7.1}
	\parab v(r,t) =\frac{\mu}{2} \sigma(t) \Big(\xi^2(r,t)+\psi(\xi(r,t)) \Big)^{-\frac{\mu}{2}-1}
	 \opa \psi((\xi(r,t))
  \ee
  for all $(r,t) \in S:=\{(\rho,\tau) \in (0,\infty)^2 \ | \ \xi(\rho,\tau) \in (\xi_0,\xi_1)\}$, where
  \be{opa}
	\opa \psi(\xi) := \Big(\xi^2+\psi \Big) \Big( \psi_{\xi\xi} + \frac{n-1}{\xi} \psi_\xi \Big)
	+ 2\kappa \psi - (\mu+\kappa) \xi \psi_\xi - \frac{\mu}{2} \psi_\xi^2
  \ee
  for $\xi \in (\xi_0,\xi_1)$.
\end{lem}
\proof
  Using $\xi_t=\frac{1}{\mu} \sigma^{\frac{1}{\mu}-1} \sigma_t  r = \frac{1}{\mu} \frac{\sigma_t}{\sigma} \xi$
  and $\sigma_t=\mu\kappa \sigma$, we compute
  \bea{7.12}
	v_t &=& \sigma_t \Big(\xi^2+\psi(\xi) \Big)^{-\frac{\mu}{2}}
	-\frac{\mu}{2} \sigma \Big(\xi^2+\psi(\xi)) \Big)^{-\frac{\mu}{2}-1}  \Big(2\xi+\psi_\xi(\xi)\Big) \xi_t
		\nn\\[1mm]
	&=& \frac{\mu}{2} \sigma \Big(\xi^2+\psi(\xi)) \Big)^{-\frac{\mu}{2}-1}
	\bigg\{ \frac{2}{\mu}  \frac{\sigma_t}{\sigma}  \Big( \xi^2+\psi(\xi) \Big)
	-\frac{1}{\mu}  \frac{\sigma_t}{\sigma}  \Big( 2\xi+\psi_\xi(\xi)\Big)  \xi \bigg\} \nn\\[1mm]
	&=& \frac{\mu}{2} \sigma \Big(\xi^2+\psi(\xi)) \Big)^{-\frac{\mu}{2}-1}
	\Big\{ 2\kappa \psi(\xi) - \kappa \xi \psi_\xi(\xi)\Big\}
  \eea
  whenever $(r,t)\in S$.
  Since $\xi_r=\sigma^\frac{1}{\mu}$, we moreover have
  $$
	(v^m)_r = \bigg( \sigma^m \Big( \xi^2 + \psi(\xi) \Big)^{-\frac{\mu m}{2}} \bigg)_r
	= - \frac{\mu m}{2} \sigma^{m+\frac{1}{\mu}} \Big(\xi^2+\psi(\xi)) \Big)^{-\frac{\mu m}{2}-1}
	\Big( 2\xi+\psi_\xi(\xi) \Big)
  $$
  as well as
  \bas
	(v^m)_{rr} &=& - \frac{\mu m}{2} \sigma^{m+\frac{2}{\mu}} \Big(\xi^2+\psi(\xi)) \Big)^{-\frac{\mu m}{2}-1}
		\Big( 2 + \psi_{\xi\xi}(\xi) \Big) \\
	& & + \frac{\mu m}{2} \Big( \frac{\mu m}{2} +1 \Big) \sigma^{m+\frac{2}{\mu}}
		\Big(\xi^2+\psi(\xi)) \Big)^{-\frac{\mu m}{2}-2} \Big( 2\xi+\psi_\xi(\xi) \Big)^2
  \eas
  at such points. Thus, in view of the identity $m+\frac{2}{\mu}=1$ and, equivalently, $\frac{\mu m}{2}=\frac{\mu}{2}-1$,
  we find that
  \bea{7.13}
	\frac{1}{m} \Big( (v^m)_{rr} + \frac{n-1}{r} (v^m)_r \Big)
	&=& -\frac{\mu}{2} \sigma \Big(\xi^2+\psi(\xi) \Big)^{-\frac{\mu}{2}} \Big( 2+ \psi_{\xi\xi}(\xi) \Big) \nn\\
	& & + \Big(\frac{\mu}{2}\Big)^2 \sigma \Big(\xi^2+\psi(\xi) \Big)^{-\frac{\mu}{2}-1}
	\Big( 2\xi + \psi_\xi(\xi) \Big)^2 \nn\\
	& &  -\frac{\mu}{2} \frac{n-1}{r} \sigma^{1-\frac{1}{\mu}}
	\Big(\xi^2+\psi(\xi) \Big)^{-\frac{\mu}{2}} \Big( 2\xi+\psi_\xi(\xi) \Big) \nn\\[2mm]
	& & \hspace*{-20mm}
	= \, \frac{\mu}{2} \sigma \Big(\xi^2+\psi(\xi) \Big)^{-\frac{\mu}{2}-1} \bigg\{
	-\Big(\xi^2+\psi(\xi) \Big) \Big( \psi_{\xi\xi} + \frac{n-1}{\xi} \psi_\xi(\xi) \Big) \nn\\
	& &
	+ \frac{\mu}{2} \Big( 2\xi+\psi_\xi(\xi)\Big)^2
	-\frac{n-1}{\xi} \Big( \xi^2+\psi(\xi) \Big) \Big( 2\xi+\psi_\xi(\xi)\Big) \bigg\} \nn\\[2mm]
	& & \hspace*{-20mm}
	= \, \frac{\mu}{2} \sigma \Big(\xi^2+\psi(\xi) \Big)^{-\frac{\mu}{2}-1} \bigg\{
	-\Big( \xi^2+\psi(\xi)\Big) \Big( \psi_{\xi\xi}(\xi) + \frac{n-1}{\xi} \psi_\xi(\xi) \Big) \nn\\
	& &
	-2(n-\mu) \xi^2 -2n\psi(\xi)
	+ 2\mu \xi \psi_\xi(\xi) + \frac{\mu}{2} \psi_\xi^2(\xi) \bigg\}
  \eea
  if $(r,t)\in S$.
  As $r\xi_r=\xi$, we finally have
  \bas
	\mu r v_r &=& \mu r \bigg\{ - \frac{\mu}{2} \sigma \Big(\xi^2+\psi(\xi) \Big)^{-\frac{\mu}{2}-1}
		 \Big( 2\xi+\psi_\xi(\xi) \Big)  \xi_r \bigg\} \\[2mm]
	&=& \frac{\mu}{2} \sigma \Big(\xi^2+\psi(\xi) \Big)^{-\frac{\mu}{2}-1}
		\Big\{ -2\mu \xi^2 - \mu \xi \psi_\xi(\xi) \Big\}
  \eas
  and therefore obtain from (\ref{7.12}) and (\ref{7.13}) that
  \bas
	\parab v &=& \frac{\mu}{2} \sigma \Big(\xi^2+\psi(\xi) \Big)^{-\frac{\mu}{2}-1}  \bigg\{
	2\kappa \psi(\xi) - \kappa \xi \psi_\xi(\xi) \\
	& & \hspace*{5mm}
	+\Big(\xi^2+\psi(\xi)\Big)  \Big( \psi_{\xi\xi}(\xi) + \frac{n-1}{\xi} \psi_\xi(\xi)\Big)
	+2(n-\mu)\xi^2 \\
	& & \hspace*{5mm}
	+2n\psi(\xi) - 2\mu \xi\psi_\xi(\xi) - \frac{\mu}{2}\psi_\xi^2(\xi)
	+2\mu \xi^2 + \mu \xi \psi_\xi(\xi)
	-2n\xi^2 - 2n\psi(\xi) \bigg\},
  \eas
  which after a straightforward rearrangement yields (\ref{7.1}). \qed
\begin{lem}\label{lem2}
  Let $n\ge 1, m>0$ and $\eps>0$. Then
  \bas
	\psi(\xi):=1-\eps\xi^2, \qquad \xi \ge 0,
  \eas
  satisfies
  \be{2.1}
	\opa \psi(\xi)=2(\kappa-n\eps) - 2(n-\mu) \eps(1-\eps) \xi^2
	\qquad \mbox{for all } \xi>0,
  \ee
  where $\opa$ is defined by {\rm (\ref{opa})}.
\end{lem}
\proof
  We directly compute
  \bas
	\opa \psi (\xi) \!\! &=& \!\!(\xi^2+1+\eps\xi^2) ( -2\eps -(n-1)  2\eps) + 2\kappa (1-\eps \xi^2)
	+ (\mu+\kappa) 2\eps \xi^2 - \frac{\mu}{2}  (-2\eps\xi)^2 \\
	&=& \!\! \Big( (1-\eps) \xi^2 +1\Big)  (-2n\eps) + 2\kappa - 2\kappa \eps \xi^2
	+2\kappa\eps\xi^2 + 2\mu\eps\xi^2 -2\mu\eps^2 \xi^2 \\
	&=& \!\! -2n\eps(1-\eps)\xi^2 -2n\eps +2\kappa +2\mu\eps\xi^2 - 2\mu\eps^2 \xi^2
	\qquad \mbox{for } \xi>0,
  \eas
  and thereby immediately obtain (\ref{2.1}). \qed

\mysection{Lower bound}\label{sect_lower}

Once we have the preparatory material, we proceed next to establish 
the lower bound for the solutions mentioned in Theorem~\ref{thm1}. 
This is the content of Proposition~\ref{lem5}. Section~\ref{sect_upper} 
will contain the proof of the corresponding upper bound, Proposition~\ref{lem6}.

\begin{lem}\label{lem4}
  Suppose that condition {\rm (\ref{mn})} on $n$ and $m$ holds, and that $l\in (\mu+2,n)$ and $\alpha \in [\alpha_-,\alpha_+]$
  with $\alpha_\pm$ given by {\rm (\ref{alpha})} and $\kappa$ as in {\rm (\ref{kappa})}.
  Then, there exist $a>0,$ $\xi_0>0$ and a positive $\upsi \in C^0([0,\infty)) \cap C^2([0,\infty) \setminus \{\xi_0\})$
  satisfying
  \be{4.1}
	\opa \upsi \le 0 \qquad \mbox{in } (0,\infty) \setminus \{\xi_0\}
  \ee
  and
  \be{4.2}
	\liminf_{\xi\nearrow \xi_0} \upsi_\xi(\xi) > \limsup_{\xi\searrow \xi_0} \upsi_\xi(\xi),
  \ee
  and such that
  \be{4.3}
	\xi^\alpha \upsi(\xi) \to a \qquad \mbox{as } \xi\to\infty.
  \ee
  Here, $\opa$ is the operator defined in {\rm (\ref{opa})}.
\end{lem}
\proof
  We let $\eps:=\frac{\kappa}{n}$ and fix $a>0$ small such that
  \be{4.4}
	\bigg\{ \Big( \frac{2\eps}{\alpha}\Big)^\frac{\alpha}{\alpha+2}
	+ \eps  \Big( \frac{\alpha}{2\eps}\Big)^\frac{2}{\alpha+2} \bigg\}  a^\frac{2}{\alpha+2} \, < \, 1.
  \ee
  Then the function $\varphi$ defined by
  \bas
	\varphi(\xi):=a\xi^{-\alpha} - 1 + \eps \xi^2, \qquad \xi>0,
  \eas
  has a unique local minimum at the point at which
  $\varphi_\xi(\xi)=-\alpha a \xi^{-\alpha-1}+2\eps\xi=0$, that is, at the point
  \bas
	\xi_{min}=\Big( \frac{\alpha a}{2\eps} \Big)^\frac{1}{\alpha+2},
  \eas
  with corresponding minimum value
  $$
	\varphi(\xi_{min}) = a  \Big( \frac{\alpha a}{2\eps} \Big)^{-\frac{\alpha}{\alpha+2}} -1
	+\eps  \Big( \frac{\alpha a}{2\eps} \Big)^\frac{2}{\alpha+2}
	= \bigg\{ \Big( \frac{2\eps}{\alpha}\Big)^\frac{\alpha}{\alpha+2}
	+ \eps  \Big( \frac{\alpha}{2\eps} \Big)^\frac{2}{\alpha+2} \bigg\} a^\frac{2}{\alpha+2} -1
	< 0
  $$
  by (\ref{4.4}). Therefore,
  \bas
	\xi_0:=\inf \{\xi>0 \ | \ \varphi(\xi) \le 0 \}
  \eas
  lies in $(0,\xi_{min})$ and we have
  \be{4.5}
	\varphi_\xi(\xi_0)<0.
  \ee
  Accordingly,
  \bas
	\upsi(\xi):=\left\{ \begin{array}{ll}
	1-\eps \xi^2 \quad & \mbox{if } \xi \in [0,\xi_0], \\[1mm]
	a\xi^{-\alpha} & \mbox{if } \xi>\xi_0,
	\end{array} \right.
  \eas
  defines a positive continuous function on $[0,\infty)$ which, by (\ref{4.5}), satisfies (\ref{4.2}) and clearly also
  fulfils (\ref{4.3}).
  Moreover, using $l>\mu+2$ and (\ref{kappa}) we find
  \bas
	(l-\mu)(n-\kappa) &=& (l-\mu)n - (l-\mu-2)(n-l)
	= l^2 - (\mu+2)l + 2n \\
	&=& \Big(l-\frac{\mu+2}{2}\Big)^2 - \frac{(\mu+2)^2}{4} + 2n
	> \Big( \mu+2-\frac{\mu+2}{2}\Big)^2 - \frac{(\mu+2)^2}{4} +2n \\
	&=& 2n>0,
  \eas
  hence $\eps \equiv \frac{\kappa}{n}<1$. Thus, recalling Lemma~\ref{lem2} and the fact that $\mu<n$, we obtain
  \be{4.6}
	\opa \upsi(\xi) = - 2(n-\mu)\eps(1-\eps)\xi^2 < 0
	\qquad \mbox{for all } \xi\in (0,\xi_0).
  \ee

  As to large $\xi$, we compute
  \bas
	\opa \upsi(\xi) &=& \Big(\xi^2 + a\xi^{-\alpha} \Big)
	\Big( \alpha(\alpha+1) a \xi^{-\alpha-2} - (n-1)\alpha a \xi^{-\alpha-2} \Big) \\
	& & + 2\kappa a \xi^{-\alpha} - (\mu+\kappa)\alpha a \xi^{-\alpha}
	-\frac{\mu}{2} \alpha^2 a^2 \xi^{-2\alpha-2} \\[2mm]
	&=& \Big\{ \alpha(\alpha+1) - (n-1)\alpha + 2\kappa - (\mu+\kappa)\alpha \Big\} a \xi^{-\alpha} \\
	& & + \Big\{\alpha(\alpha+1) - (n-1)\alpha - \frac{\mu}{2} \alpha^2 \Big\}  a^2 \xi^{-2\alpha-2} \\[2mm]
	&=& \Big\{ \alpha^2-(n-2-\mu-\kappa) \alpha + 2\kappa \Big\}  a\xi^{-\alpha} \\
	& & - \Big\{ \frac{\mu-2}{2} \alpha^2 + (n-2)\alpha\Big\}  a^2 \xi^{-2\alpha-2}
	\qquad \mbox{for } \xi>\xi_0.
  \eas
  Since $\alpha^2-(n-2-\mu-\kappa)\alpha+2\kappa \le 0$ due to our assumption $\alpha \in [\alpha_-,\alpha_+]$, and
  since $\mu=\frac{2}{1-m}>2$ and $n>2$, from this we immediately infer that
  \bas
	\opa \upsi(\xi) < 0 \qquad \mbox{for all } \xi>\xi_0,
  \eas
  which combined with (\ref{4.6}) proves (\ref{4.1}).
\qed
\begin{prop}\label{lem5}
  Assume again condition {\rm (\ref{mn})}. Suppose that $v_0 \in C^0([0,\infty))$ is positive and fulfils
  \begin{equation}\label{5.1}
	v_0(r) \ge r^{-\mu}-c_0 r^{-l} \qquad \mbox{for all } r\ge 1
  \end{equation}
  with some $c_0>0$ and $l$ as in {\rm (\ref{eq:l})}.
 Let $v$ denote the solution
  of {\rm (\ref{0})}. Then:\\
  (i) \ There exists $C_1>0$ such that
  \begin{equation}\label{5.2}
	v(0,t) \ge C_1 \, e^{\frac{\mu(l-\mu-2)(n-l)}{l-\mu} t} \qquad \mbox{for all } t\ge 0.
  \end{equation}
  (ii) \ For each $r_0>$ one can find $C_2>0$ fulfilling
  \begin{equation}\label{200}
	v(r,t) \ge r^{-\mu}- C_2 \, e^{-(l-\mu-2)(n-l)t} \qquad \mbox{for all $r\ge r_0$ and } t \ge 0.
  \end{equation}
\end{prop}
\proof
  (i) \ As $l$ is as in (\ref{eq:l}),
  we can apply Lemma~\ref{lem4} to $\alpha_-=l-\mu-2$ and obtain $\xi_0>0$
  and a function $\upsi$ with the properties provided by that lemma.
  Since $v_0$ is positive, there exists $c_1>0$ such that
  \be{5.3}
	v_0(r) \ge c_1 \qquad \mbox{for all } r\in [0,1].
  \ee
  Let us pick $c_2>0$ large such that with $c_0$ as in (\ref{5.1}) we have
  \be{5.4}
	(1+c_2)^{-\frac{\mu}{2}} \le 1-c_0,
  \ee
  and then choose $\hat\xi>0$ small fulfilling
  \be{5.5}
	\xi^{-2} \upsi(\xi) \ge c_2 \qquad \mbox{for all } \xi \in (0,\hat\xi),
  \ee
  which is possible since $\upsi$ is positive on $[0,\infty)$.
  We next define
  \be{5.6}
	z_0:=\max_{\xi \in [\hat\xi,\infty)} \xi^{-2} \upsi(\xi)
  \ee
  and let $c_3>0, c_4>0$ and $c_5>0$ be small enough such that
  \be{5.7}
	(1+z)^{-\frac{\mu}{2}} \le 1-c_3 z \qquad \mbox{for all } z\in [0,z_0]
  \ee
  and
  \be{5.8}
	\xi^2 + \upsi(\xi) \ge c_4 \qquad \mbox{for all } \xi \ge 0
  \ee
  as well as
  \be{5.9}
	\upsi(\xi) \ge c_5 \xi^{-\alpha} \qquad \mbox{for all } \xi \ge \hat\xi,
  \ee
  where we make use of (\ref{4.3}) and, again, the positivity of $\upsi$.
  Finally, we take a small number $\sigma_0>0$ satisfying
  \be{5.10}
	\sigma_0 \le c_1 c_4^\frac{\mu}{2}
  \ee
  and
  \be{5.11}
	\sigma_0 \le \Big(\frac{c_3 c_5}{c_0}\Big)^\frac{\mu}{l-\mu}
  \ee
  and define
  \bas
	\sigma(t) := \sigma_0  e^{\mu\kappa t}, \qquad t\ge 0,
  \eas
  with $\kappa=\frac{(l-\mu-2)(n-l)}{l-\mu}$ as in (\ref{kappa}).
  Then
  \bas
	\uv(r,t):=\sigma(t)
\Big(\xi^2(r,t)+\upsi(\xi(r,t))\Big)^{-\frac{\mu}{2}}, \qquad r\ge 0, \ t\ge 0,
  \eas
  with $\xi(r,t):=\sigma^\frac{1}{\mu}(t)r$, is continuous in $[0,\infty)^2$ and smooth at each point
  $(r,t) \in [0,\infty)^2$ where $r \ne r_0(t):=\xi_0 \sigma^{-\frac{1}{\mu}}(t)$.
  An application of Lemma~\ref{lem7} and Lemma~\ref{lem4} shows that with $\opa$ as in (\ref{opa}),
  \bas
	\parab \uv = \frac{\mu}{2} \sigma(t)  \Big(\xi^2(r,t)+\upsi(\xi(r,t))\Big)^{-\frac{\mu}{2}-1}
	 \opa \upsi (\xi(r,t)) \, \le 0 \qquad \mbox{whenever } r\ne r_0(t),
  \eas
  which implies that $\uv$ is a subsolution of (\ref{0}) in the Nagumo sense, because from (\ref{4.2}) we infer that
  \bas
	\limsup_{r\nearrow r_0(t)} \uv_r(r,t) < \liminf_{r\searrow r_0(t)} \uv_r(r,t) \qquad \mbox{for all } t\ge 0.
  \eas
  Accordingly, if we can show that $\uv$ does not exceed $v$ initially then the comparison principle will tell us that
  \be{201}
	v \ge \uv \qquad \mbox{in } [0,\infty)^2
  \ee
  and, in particular,
  \bas
	v(0,t) \ge \uv(0,t) = \sigma(t) =\sigma_0  e^{\mu\kappa t} \qquad \mbox{for all } t\ge 0,
  \eas
  which will yield (\ref{5.2}) in view of the definition of $\kappa$.
  It thus remains to show that
  \be{5.14}
	\uv(r,0) \le v_0(r) \qquad \mbox{for all } r \ge 0.
  \ee
  To this end, we first consider the case $r \le 1$, when (\ref{5.8}),
   (\ref{5.10}) and (\ref{5.3}) imply
  \bea{5.15}
	\uv(r,0) &=& \sigma_0  \Big( \xi^2(r,0)+\upsi(\xi(r,0))\Big)^{-\frac{\mu}{2}} \nn\\
	&\le& \sigma_0  c_4^{-\frac{\mu}{2}}
	\le c_1
	\le v_0(r) \qquad \mbox{for all } r\in [0,1].
  \eea
  Next, if $r\ge 1$ is such that $r \ge \hat\xi  \sigma_0^{-\frac{1}{\mu}}$ then $\xi(r,0)=\sigma_0^\frac{1}{\mu}r
  \ge \hat\xi$, and hence from (\ref{5.6}), (\ref{5.7}) and (\ref{5.9}) we obtain
  \bas
	\uv(r,0) &=& r^{-\mu}  \Big( 1+ \xi^{-2}(r,0) \upsi(\xi(r,0)) \Big)^{-\frac{\mu}{2}}
	\le r^{-\mu}  \Big( 1-c_3 \xi^{-2}(r,0) \upsi(\xi(r,0)) \Big) \\
	&\le& r^{-\mu} \Big( 1-c_3 c_5 \xi^{-\alpha-2}(r,0)\Big)
	= r^{-\mu} - c_3 c_5 \sigma_0^{-\frac{\alpha+2}{\mu}} r^{-\mu-\alpha-2},
  \eas
  which in view of our choice $\alpha=l-\mu-2$, (\ref{5.11}) and (\ref{5.1}) gives
  \bea{5.16}
	\uv(r,0) &\le& r^{-\mu} - c_3 c_5 \sigma_0^{-\frac{l-\mu}{\mu}} r^{-l}
	\le r^{-\mu} - c_0 r^{-l} \nn\\
	&\le& v_0(r)
	\qquad \mbox{for all } r\ge \min\Big\{1,\hat \xi \sigma_0^{-\frac{1}{\mu}}\Big\}.
  \eea
  Finally, if $r\ge 1$ is such that $r \le \hat \xi \sigma_0^{-\frac{1}{\mu}}$ then we use (\ref{5.5}) and (\ref{5.4})
  to see that
  \bas
	\uv(r,0) &=& r^{-\mu}  \Big( 1+\xi^{-2}(r,0) \upsi(\xi(r,0)) \Big)^{-\frac{\mu}{2}}
	\le r^{-\mu}  (1+c_2)^{-\frac{\mu}{2}} \\
	&\le& r^{-\mu}  (1-c_0) \qquad \mbox{for all } r \le \hat \xi \sigma_0^{-\frac{1}{\mu}},
  \eas
  whereas by (\ref{5.1}),
  $$
	v_0(r) \ge r^{-\mu}  \Big(1-c_0 r^{-(l-\mu)} \Big)
	\ge r^{-\mu} (1-c_0) \qquad \mbox{for all } r \ge 1,
  $$
  so that $\uv(r,0) \le v_0(r)$ also holds if $1 \le r \le \hat \xi \sigma_0^{-\frac{1}{\mu}}$. Together with (\ref{5.15})
  and (\ref{5.16}) this establishes (\ref{5.14}).\abs
  (ii) \ To see (\ref{200}), we observe that in view of (\ref{4.3}) there exists $c_6>0$ such that
  \bas
	\upsi(\xi) \le c_6 \xi^{-\alpha} \qquad \mbox{for all } \xi>0,
  \eas
  where still $\alpha=\alpha_-=l-\mu-2$. Then by (\ref{201}) and the convexity of $z\mapsto (1+z)^{-\frac{\mu}{2}}$ for
  $z\ge 0$ we obtain
  \bas
	v(r,t) &\ge& \uv(r,t)
	= r^{-\mu} \Big(1+\xi^{-2} \upsi(\xi)\Big)^{-\frac{\mu}{2}}
	\ge r^{-\mu} - \frac{\mu}{2} r^{-\mu} \xi^{-2} \upsi(\xi) \\
	&\ge& r^{-\mu} - \frac{\mu}{2} c_6 r^{-\mu} \xi^{-\alpha-2}
	= r^{-\mu}-\frac{\mu}{2} c_6 \sigma^{-\frac{\alpha+2}{\mu}} r^{-\mu-\alpha-2}\\
	&=& r^{-\mu} - \frac{\mu}{2} c_6 \sigma_0^{-\frac{l-\mu}{\mu}}
        r^{-l} \, e^{-(l-\mu-2)(n-l)t} \qquad \mbox{for all } r>0 \mbox{ and }
	t>0.
  \eas
  Given $r_0>0$, this easily yields (\ref{200}) upon an obvious choice of $C_2$.
\qed
\mysection{Upper bound and proof of Theorem~\ref{thm1}}\label{sect_upper}
\begin{lem}\label{lem1}
  Assume {\rm (\ref{mn})}, and let $l\in (\mu+2,n)$ and $\alpha_-$ be as defined in
  {\rm (\ref{alpha})}
  with $\kappa$ given by {\rm (\ref{kappa})}.
  Then there exist $\beta>\alpha_-$
  and $C_\beta>0$ with the following property: Suppose that $A>0$ and $B>0$ are such that
  \be{1.1}
	\frac{B^{\alpha_-+2}}{A^{\beta+2}} \ge C_\beta,
  \ee
  and let
  \be{1.4}
	\xi_1:=\Big(\frac{\beta B}{\alpha_- A}\Big)^\frac{1}{\beta-\alpha_-}.
  \ee
  Then the function $\opsio$ defined by
  \be{1.2}
	\opsio(\xi):=A\xi^{-\alpha_-} -B\xi^{-\beta}, \qquad \xi>0,
  \ee
  satisfies
  \be{1.444}
	(\opsio)_\xi(\xi_1)=0 \qquad \mbox{and} \qquad (\opsio)_\xi(\xi) <0 \quad \mbox{for all } \xi>\xi_1,
  \ee
  and moreover we have
  \be{1.3}
	\opa \opsio(\xi) \ge 0 \qquad \mbox{for all } \xi > \xi_1
  \ee
  with $\opa$ given by {\rm (\ref{opa})}.
\end{lem}
\proof
  Throughout the proof, let us abbreviate $\alpha:=\alpha_-$ for convenience.
  Then since $l\in (\mu+2,n)$, we have $0<\alpha < \alpha_+$, and
  \bas
	p(\beta):=\beta^2-(n-2-\mu-\kappa)\beta + 2\kappa, \qquad \beta \in \R,
  \eas
  has the properties
  \be{1.11}
	p(\alpha)=0 \qquad \mbox{and} \qquad p(\beta) < 0 \quad \mbox{for all } \beta \in (\alpha,\alpha_+).
  \ee
  As furthermore
  \bas
	q(\beta):=\beta(\beta+2-n)-\alpha(\alpha+2-n) + \mu\alpha\beta, \qquad \beta\in\R,
  \eas
  evidently satisfies $q(\alpha)=\mu\alpha^2>0$, by continuity we can choose some $\beta>\alpha$ such that
  \be{1.12}
	\beta \le 2\alpha+2, \qquad p(\beta)<0 \qquad \mbox{and} \qquad q(\beta)>0.
  \ee
  With this value of $\beta$ fixed, we pick $C_\beta>0$ large such that
  \be{1.13}
	\frac{2c_2}{c_1}  \Big( \frac{\alpha}{\beta} \Big)^\frac{2\alpha+2-\beta}{\beta-\alpha}
	 \Big( \frac{1}{C_\beta}\Big)^\frac{1}{\beta-\alpha} \, \le 1
  \ee
  and
  \be{1.14}	
	\frac{2c_3}{c_1}  \Big( \frac{\alpha}{\beta}\Big)^\frac{\beta+2}{\beta-\alpha}
	 \Big(\frac{1}{C_\beta}\Big)^\frac{1}{\beta-\alpha} \, \le 1,
  \ee
  where
  \be{1.144}
	c_1:=-p(\beta), \qquad c_2:=\frac{\mu-2}{2}\alpha^2 + (n-2)\alpha \qquad \mbox{and} \qquad
	c_3:=\frac{\mu-2}{2} \beta^2 + (n-2)\beta
  \ee
  are all positive according to (\ref{1.12}) and the inequalities $\mu=\frac{2}{1-m}>2$ and $n>2$.
  Then, given $A>0$ and $B>0$ fulfilling (\ref{1.1}), we let $\opsio$ be defined by (\ref{1.2}) and compute
  \bea{1.44}
	(\opsio)_\xi(\xi) &=& -\alpha A \xi^{-\alpha-1} + \beta B \xi^{-\beta-1}
	\qquad \mbox{and} \nn\\
	(\opsio)_{\xi\xi}(\xi) &=& \alpha(\alpha+1) A\xi^{-\alpha-2} - \beta(\beta+1) B\xi^{-\beta-2}
  \eea
  for $\xi>0$.
  From this it can easily be deduced that in fact $\opsio$ attains its maximum at $\xi=\xi_1$ and decreases
  on $(\xi_1,\infty)$, where $\xi_1$ is as in (\ref{1.4}). Using (\ref{1.44}) we furthermore obtain
  \bas
	\opa \opsio(\xi) &=& \big( \xi^2 + A\xi^{-\alpha} - B\xi^{-\beta} \big)
	\Big\{ \alpha(\alpha+1) A \xi^{-\alpha-2} - \beta(\beta+1) B \xi^{-\beta-2} \\
	& & \hspace*{40mm} -(n-1)\alpha A \xi^{-\alpha-2} + (n-1)\beta B \xi^{-\beta-2} \Big\} \\
	& & + 2\kappa A \xi^{-\alpha} - 2\kappa B \xi^{-\beta}
	 - (\mu+\kappa) \alpha A \xi^{-\alpha} + (\mu+\kappa)\beta B \xi^{-\beta} \\
	& & - \frac{\mu}{2} \Big(\alpha^2 A^2 \xi^{-\alpha-2} - 2\alpha \beta AB \xi^{-\alpha-\beta-2}
	+\beta^2 B^2 \xi^{-2\beta-2}\Big) \\[2mm]
	&=& \Big\{ \alpha(\alpha+1) - (n-1)\alpha + 2\kappa -(\mu+\kappa)\alpha \Big\}  A \xi^{-\alpha} \\
	& & + \Big\{ -\beta(\beta+1)+(n-1)\beta) - 2\kappa + (\mu+\kappa)\beta \Big\}  B\xi^{-\beta} \\
	& & + \Big\{ \alpha(\alpha+1) - (n-1)\alpha-\frac{\mu}{2} \alpha^2 \Big\}  A^2 \xi^{-2\alpha-2} \\
	& & + \Big\{ \beta(\beta+1) - (n-1)\beta - \frac{\mu}{2} \beta^2 \Big\}  B^2 \xi^{-2\beta-2} \\
	& & + \Big\{ -\alpha(\alpha+1) + (n-1)\alpha + \beta(\beta+1) - (n-1)\beta + \mu\alpha\beta \Big\}
	 AB \xi^{-\alpha-\beta-2} \\[2mm]
	&=& p(\alpha)  A\xi^{-\alpha} - p(\beta)  B\xi^{-\beta}
	 - \Big\{\frac{\mu-2}{2} \alpha^2 + (n-2)\alpha \Big\}  A^2\xi^{-2\alpha-2}\\
	& & - \Big\{ \frac{\mu-2}{2} \beta^2 + (n-2) \beta \Big\}  B^2 \xi^{-2\beta-2}
	 + q(\beta)  AB \xi^{-\alpha-\beta-2}
	\qquad \mbox{for all } \xi>0.
  \eas
  Now (\ref{1.11}) and (\ref{1.12}) imply that the first term on the right vanishes and that the last is nonnegative,
  because $A$ and $B$ are positive. Hence, recalling (\ref{1.144}) we arrive at the inequality
  \be{1.5}
	\opa \opsio(\xi) \ge c_1 B \xi^{-\beta} - c_2 A^2 \xi^{-2\alpha-2} - c_3 B^2 \xi^{-2\beta-2}
	\qquad \mbox{for all } \xi>0.
  \ee
  Now if $\xi\ge \xi_0$ then (\ref{1.13}) along with (\ref{1.1}) and our restriction $\beta \le 2\alpha+2$ ensures that
  \bea{1.6}
	\frac{c_2 A^2 \xi^{-2\alpha-2}}{\frac{1}{2} c_1 B \xi^{-\beta}}
	&=& \frac{2c_2}{c_1}  \frac{A^2}{B}  \xi^{\beta-2\alpha-2}
	\le \frac{2c_2}{c_1}  \frac{A^2}{B}  \Big( \frac{\beta B}{\alpha A}
		\Big)^\frac{\beta-2\alpha-2}{\beta-\alpha} \nn\\
	&=& \frac{2c_2}{c_1}  \Big(\frac{\alpha}{\beta}\Big)^\frac{2\alpha+2-\beta}{\beta-\alpha}
	\Big(\frac{A^{\beta+2}}{B^{\alpha+2}} \Big)^\frac{1}{\beta-\alpha}
	\le 1 \qquad \mbox{for all } \xi\ge \xi_1.
  \eea
  Moreover, for such $\xi$ we find
  \bas
	\frac{c_3 B^2 \xi^{-2\beta-2}}{\frac{1}{2} c_1 B \xi^{-\beta}}
	&=& \frac{2c_3}{c_1}  B  \xi^{-\beta-2}
	\le \frac{2c_3}{c_1}  B  \Big( \frac{\beta B}{\alpha A} \Big)^\frac{-\beta-2}{\beta-\alpha} \\
	&=& \frac{2c_3}{c_1}  \Big( \frac{\alpha}{\beta} \Big)^\frac{\beta+2}{\beta-\alpha}
	 \Big(\frac{A^{\beta+2}}{B^{\alpha+2}} \Big)^\frac{1}{\beta-\alpha}
	\le 1 \qquad \mbox{for all } \xi \ge \xi_1
  \eas
  by (\ref{1.14}). Together with (\ref{1.6}) and (\ref{1.5}), this shows that indeed $\opa \opsio \ge 0$ for all
  $\xi\ge \xi_1$, as claimed.
\qed
\begin{lem}\label{lem3}
  Suppose that {\rm (\ref{mn})} holds. Let $l \in (\mu+2,n)$ and $\alpha_-$ be as in
  {\rm (\ref{alpha})} with $\kappa$ given by {\rm (\ref{kappa})}.
  Then there exist $A>0$, $\xi_1>0$ and a positive function $\opsi \in C^0([0,\infty)) \cap
  C^2([0,\infty) \setminus \{\xi_1\})$ such that with $\opa$ as in {\rm (\ref{opa})},
  \be{3.1}
	\opa \opsi \ge 0 \qquad \mbox{in } (0,\infty) \setminus \{\xi_1\}
  \ee
  and
  \be{3.2}
	\limsup_{\xi\nearrow \xi_1} \opsi_\xi(\xi) <
          \liminf_{\xi \searrow \xi_1} \opsi_\xi(\xi)
  \ee
  as well as
  \be{3.2222}
	\xi^{\alpha_-} \opsi(\xi) \to A \qquad \mbox{as } \xi\to\infty.
  \ee
\end{lem}
\proof
  Again we write $\alpha:=\alpha_-$ for simplicity.
  Since $\mu<n$, it is possible to fix $c_1 \in (0,1)$ such that
  \be{3.222}
	c_1 \le \sqrt{\frac{\kappa}{2(n-\mu)}},
  \ee
  and since $l\in (\mu+2,n)$, there exist $\beta>\alpha$ and $C_\beta>0$ such that the conclusion of
  Lemma~\ref{lem1}
  holds. We now define
  \be{3.22}
	c_2:=\Big(\frac{\alpha}{\beta}\Big)^\frac{\alpha}{\beta-\alpha}
	- \Big(\frac{\alpha}{\beta}\Big)^\frac{\beta}{\beta-\alpha},
  \ee
  which is positive because $0<\alpha<\beta$, and
  \be{3.3}
	K:=\Big(\frac{1-c_1^2}{c_2}\Big)^{\beta-\alpha}.
  \ee
  Next, we pick $A>0$ large fulfilling
  \be{3.4}
	A \ge \Big\{ \frac{2n}{\kappa}  c_1^2  \Big(\frac{\alpha}{\beta}\Big)^\frac{2}{\beta-\alpha}
	 K^\frac{2}{\alpha(\beta-\alpha)} \Big\}^\frac{\alpha}{2}
  \ee
  and
  \be{3.5}
	A \ge \Big( C_\beta  K^\frac{\alpha+2}{\alpha} \Big)^\frac{\alpha}{2(\beta-\alpha)}
  \ee
  and let
  \be{3.6}
	B:=\Big(\frac{A^\beta}{K}\Big)^\frac{1}{\alpha}
  \ee
  and
  \be{3.7}
	\eps:= c_1^2  \Big( \frac{\alpha A}{\beta B} \Big)^\frac{2}{\beta-\alpha}
  \ee
  as well as
  \be{3.8}
	\xi_1:= \Big( \frac{\beta B}{\alpha A} \Big)^\frac{1}{\beta-\alpha},
  \ee
  so that
  \be{3.9}
	\eps \xi_1^2 = c_1^2 \Big(\frac{\alpha A}{\beta B} \Big)^\frac{2}{\beta-\alpha}
	 \Big(\frac{\beta B}{\alpha A} \Big)^\frac{2}{\beta-\alpha} \, = c_1^2.
  \ee
  Then the function $\opsi:[0,\infty) \to \R$ given by
  \bas
	\opsi(\xi) := \left\{ \begin{array}{ll}
	\opsii(\xi):=1-\eps \xi^2 & \mbox{if } \xi \in [0,\xi_1], \\[1mm]
	\opsio(\xi)=A\xi^{-\alpha}-B\xi^{-\beta} \quad & \mbox{if } \xi \in (\xi_1,\infty),
	\end{array} \right.
  \eas
  is continuous on $[0,\infty)$, because (\ref{3.9}) ensures that
  \bas
	\opsii(\xi_1)=1-c_1^2,
  \eas
  whereas invoking (\ref{3.22}), (\ref{3.6}) and (\ref{3.3}) we find
  $$
	\opsio(\xi_1) = A  \Big(\frac{\beta B}{\alpha A} \Big)^{-\frac{\alpha}{\beta-\alpha}}
	-B  \Big(\frac{\beta B}{\alpha A} \Big)^{-\frac{\beta}{\beta-\alpha}}
	= c_2 A^\frac{\beta}{\beta-\alpha} B^{-\frac{\alpha}{\beta-\alpha}}
	= c_2 K^\frac{1}{\beta-\alpha}
	= 1-c_1^2.
  $$
  Note that since $c_1<1$, this also implies that $\opsi$ is positive on $[0,\infty)$.
  Next, from Lemma~\ref{lem2} and (\ref{3.9}) we obtain
  \bas
	\opa \opsii(\xi) &=& 2(\kappa-n\eps) - 2(n-\mu)\eps(1-\eps) \xi^2 \\
	&\ge& 2(\kappa-n\eps) - 2(n-\mu)\eps \xi_1^2 \\
	&=& 2(\kappa-n\eps) - 2(n-\mu) c_1^2
	\qquad \mbox{for all } \xi \in (0,\xi_1),
  \eas
  where by (\ref{3.7}), (\ref{3.6}) and (\ref{3.4}),
  $$
	\eps = c_1^2  \Bigg( \frac{\alpha A}{\beta  \big(\frac{A^\beta}{K}\big)^\frac{1}{\alpha}}
		\Bigg)^\frac{2}{\beta-\alpha}
	= c_1^2  \Big(\frac{\alpha}{\beta}\Big)^\frac{2}{\beta-\alpha}  K^\frac{2}{\alpha(\beta-\alpha)}
	 A^{-\frac{2}{\alpha}}
	\le \frac{\kappa}{2n} \, .
  $$
  Hence,
  \be{3.10}
	\opa \opsii(\xi) \ge 2\Big(\kappa-\frac{\kappa}{2}\Big) - 2(n-\mu) c_1^2
	\ge 0
	\qquad \mbox{for all } \xi \in (0,\xi_1)
  \ee
  according to (\ref{3.222}).\\
  Now the requirement (\ref{3.5}) along with (\ref{3.6}) guarantees that
  \bas
	\frac{B^{\alpha+2}}{A^{\beta+2}} = K^{-\frac{\alpha+2}{\alpha}} A^\frac{2(\beta-\alpha)}{\alpha} \ge C_\beta,
  \eas
  so that Lemma~\ref{lem1} becomes applicable to tell us that
  \be{3.11}
	\opa \opsio(\xi) \ge 0 \qquad \mbox{for all } \xi>\xi_1
  \ee
  as well as
  \be{3.12}
	(\opsio)_\xi(\xi_1)=0.
  \ee
  As a consequence of (\ref{3.10}) and (\ref{3.11}), we see that (\ref{3.1}) holds, while (\ref{3.12}) combined with the fact that
  \bas
	(\opsii)_\xi(\xi_1) = - 2\eps\xi_1^2 < 0
  \eas
  yields (\ref{3.2}). The assertion (\ref{3.2222}) immediately results from the definition of $\opsi$.
\qed
%
%
%
\begin{prop}\label{lem6}
  Suppose that {\rm (\ref{mn})} holds, and that $v$ is the solution of {\rm
  (\ref{0})}, where
  the initial data $v_0 \in C^0([0,\infty))$ are nonnegative and such that
  there exist \ $l$ as in {\rm (\ref{eq:l})} and $c_1>0$ fulfilling
  \be{6.1}
	v_0(r) \le r^{-\mu} - c_1 r^{-l} \qquad \mbox{for all } r\ge 1,
  \ee
  and which in addition satisfies
  \be{429}
	v_0(r)<r^{-\mu} \qquad \mbox{for all } r>0.
  \ee
  i) \ There exists $C_1>0$ such that
  \be{6.2}
	v(r,t) \le C_1 \, e^{\frac{\mu(l-\mu-2)(n-l)}{l-\mu}  t}
	\qquad \mbox{for all $r\ge 0$ and } t\ge 0.
  \ee
  ii) \ For all $r_0>0$ there exists $C_2>0$ with the property
  \be{4311}
	v(r,t) \le r^{-\mu} - C_2 \, e^{-(l-\mu-2)(n-l)t} \qquad \mbox{for all $r\ge r_0$ and } t\ge 0.
  \ee
\end{prop}
%
\proof
  i) \
  Since $l$ is as in (\ref{eq:l}), the number $\alpha:=\alpha_-$ satisfies $\alpha=l-\mu-2$ by
  Lemma~\ref{lem_alpha}. Hence, applying Lemma~\ref{lem3} we find $\xi_1>0$ and a positive function $\opsi \in C^0([0,\infty))
  \cap C^2([0,\infty) \setminus \{\xi_1\})$ with the properties (\ref{3.1}) and (\ref{3.2}) and such that
  \be{6.3}
	\xi^\alpha \opsi(\xi) \le c_2 \qquad \mbox{for all } \xi \ge 0
  \ee
  with some $c_2>0$. Taking $c_3>0$ large such that
  \be{6.4}
	v_0(r) \le c_3 \qquad \mbox{for all } r\ge 0,
  \ee
  we can find $r_0>0$ small enough fulfilling
  \be{6.5}
	r_0 \le (2c_3)^{-\frac{1}{\mu}}
  \ee
  and then, by (\ref{6.1}) and (\ref{429}), fix $c_4 \in (0,c_1]$ such that
  \be{6.6}
	v_0(r) \le r^{-\mu} - c_4 r^{-l} \qquad \mbox{for all } r \ge r_0.
  \ee
  We pick $\hat\xi>0$ and $c_5>0$ sufficiently large satisfying
  \be{6.7}
	\hat\xi \ge (\mu c_2)^\frac{1}{\alpha+2}
  \ee
  and
  \be{6.8}
	\xi^2 + \opsi(\xi) \le c_5 \qquad \mbox{for all } \xi \in [0,\hat\xi]
  \ee
  and finally choose a large number $\sigma_0>0$ with
  \be{6.9}
	\sigma_0 \ge c_3 c_5^\frac{\mu}{2}
  \ee
  and
  \be{6.10}
	\sigma_0 \ge \Big( \frac{\mu c_2}{2c_4} \Big)^\frac{\mu}{l-\mu}.
  \ee
  We now define
  \bas
	\sigma(t):=\sigma_0  e^{\mu\kappa t}, \qquad t\ge 0,
  \eas
  and
  \bas
	\ov(r,t):=\sigma(t)  \Big( \xi^2(r,t)+\opsi(\xi(r,t)) \Big)^{-\frac{\mu}{2}},
	\qquad r\ge 0, \ t\ge 0,
  \eas
  again with $\xi(r,t):=\sigma^\frac{1}{\mu}(t) r$, and claim that
  \be{6.100}
	\ov(r,0) \ge v_0(r) \qquad \mbox{for all } r\ge 0.
  \ee
  Indeed, if $r \le \hat \xi  \sigma_0^{-\frac{1}{\mu}}$ then $\xi(r,0) \le \hat \xi$ and hence
  (\ref{6.8}), (\ref{6.9}) and (\ref{6.4}) imply that
  \be{6.11}
	\ov(r,0) = \sigma_0  \Big( \xi^2(r,0) + \opsi(\xi(r,0))\Big)^{-\frac{\mu}{2}}
	\ge \sigma_0  c_5^{-\frac{\mu}{2}}
	\ge c_3
	\ge v_0(r) \qquad \mbox{for } r \le \hat\xi  \sigma_0^{-\frac{1}{\mu}}.
  \ee
  Next, in the case when $\hat\xi  \sigma_0^{-\frac{1}{\mu}} \le r \le r_0$ we have
  $\xi(r,0) \ge \hat\xi$, so that using the convexity of $0 \le z \mapsto (1+z)^{-\frac{\mu}{2}}$ along with
  (\ref{6.3}), (\ref{6.7}), (\ref{6.5}) and (\ref{6.4}), we can estimate
  \bea{6.12}
	\ov(r,0) &=& r^{-\mu}  \Big( 1+ \xi^{-2}(r,0) \opsi(\xi(r,0)) \Big)^{-\frac{\mu}{2}}
	\ge r^{-\mu}  \Big( 1-\frac{\mu}{2} \xi^{-2}(r,0) \opsi(\xi(r,0)) \Big) \nn\\
	&\ge& r^{-\mu}  \Big( 1-\frac{\mu}{2} c_2 \xi^{-\alpha-2}(r,0)\Big)
	\ge r^{-\mu}  \Big( 1-\frac{\mu}{2} c_2  \frac{1}{\mu c_2} \Big)
	= \frac{1}{2} r^{-\mu} \nn\\
	&\ge& \frac{1}{2} r_0^{-\mu}
	\ge c_3
	\ge v_0(r) \qquad \mbox{if } \hat\xi  \sigma_0^{-\frac{1}{\mu}} \le r \le r_0.
  \eea
  Finally, for $r \ge \hat \xi  \sigma_0^{-\frac{1}{\mu}}$ fulfilling $r\ge r_0$, by the same convexity argument
  in conjunction with the fact that $\alpha=l-\mu-2$, from (\ref{6.10}) and (\ref{6.6}) we
  have
  \bas
	\ov(r,0) &\ge& r^{-\mu}  \Big( 1-\frac{\mu}{2} c_2  \xi^{-\alpha-2}(r,0) \Big)
	= r^{-\mu} - \frac{\mu}{2} c_2 \sigma_0^{-\frac{\alpha+2}{\mu}} r^{-l} \\
	&\ge& r^{-\mu} - c_4 r^{-l}
	\qquad \mbox{for } r\ge \max \Big\{ r_0, \hat\xi  \sigma_0^{-\frac{1}{\mu}} \Big\}.
  \eas
  Together with (\ref{6.11}) and (\ref{6.12}), this proves (\ref{6.100}). Since by
  Lemma~\ref{lem7} and Lemma~\ref{lem2},
  recalling (\ref{opa}) we have
  \bas
	\parab \ov = \frac{\mu}{2} \sigma(t)  \Big( \xi^2(r,t) + \opsi(\xi(r,t)) \Big)^{-\frac{\mu}{2}-1}
		\opa \opsi(r,t) \ge 0
	\quad \mbox{whenever } r \ne r_1(t):=\xi_1  \sigma^{-\frac{1}{\mu}}(t)
  \eas
  and
  \bas
	\liminf_{r\nearrow r_1(t)} \ov_r(r,t) > \limsup_{r\searrow r_1(t)} \ov_r(r,t) \qquad \mbox{for all } t\ge 0
  \eas
  according to (\ref{3.2}), the comparison principle applies to yield
  \be{99}
	\ov \ge v \qquad \mbox{for all $r\ge 0$ and } t\ge 0.
  \ee
  This immediately leads to (\ref{6.2}).\abs
  ii) \ To obtain (\ref{4311}), we fix $r_0>0$ and first pick $c_5>0$ small enough fulfilling
  \be{101}
	(1+z)^{-\frac{\mu}{2}} \le 1-c_5 z \qquad \mbox{for all } z \in [0,1],
  \ee
  and then fix $t_0>0$ large such that
  \be{100}
	I:=c_2 \sigma_0^{-\frac{l-\mu}{\mu}} r_0^{-(l-\mu)} \, e^{-(l-\mu-2)(n-l)t_0} \le 1,
  \ee
  where $c_2$ and $\sigma_0$ are as determined by (\ref{6.3}), (\ref{6.9}) and (\ref{6.10}).
  Then for all $r\ge r_0$ and $t\ge t_0$, still writing $\alpha=\alpha_-=l-\mu-2$ we have
  \bas
	c_2 \xi^{-\alpha-2}(r,t) \le c_2 \xi^{-\alpha-2}(r_0,t_0)=I \le 1
  \eas
  and thus
  \bas
	(1+c_2 \xi^{-\alpha-2})^{-\frac{\mu}{2}} \le 1-c_2c_5 \xi^{-\alpha-2}
  \eas
  by (\ref{101}).
  Therefore, (\ref{99}) and (\ref{6.3}) entail that for such $r$ and $t$ we have
  \bas
	v(r,t) &\le& \ov(r,t)
	= r^{-\mu} \Big(1+\xi^{-2} \opsi(\xi)\Big)^{-\frac{\mu}{2}} \\
	&\le& r^{-\mu} (1+c_2 \xi^{-\alpha-2})^{-\frac{\mu}{2}}
	\le r^{-\mu} (1-c_2 c_5 \xi^{-\alpha-2}) \\
	&=& r^{-\mu} - c_2 c_5  \sigma^{-\frac{\alpha+2}{\mu}}  r^{-\mu-\alpha-2}
	= r^{-\mu} - c_2 c_5 \sigma_0^{-\frac{l-\mu}{\mu}} r^{-l} \, e^{-(l-\mu-2)(n-l)t}.
  \eas
  This shows that (\ref{4311}) is valid for some sufficiently large $C_2>0$.
\qed

\

As we have said, Propositions ~\ref{lem5} and ~\ref{lem6} together imply Theorem~\ref{thm1}.

\mysection{Comments and open problems}
\label{sec.cop}

\noindent {\bf 1.} The construction of the new extinction rates for $m\in [m_*, m_c) $ is open.  The relevance of $m_*$ in the asymptotic analysis of stability of the Barenblatt solutions has been documented in
\cite{BBDGV, BGV, PNAS}.

\medskip

\noindent {\bf 2.} We have not performed the analysis of positive perturbations of the tail of the singular solution $V_0$. Preliminary calculations show that we can have in that case global grow-up if the perturbation is large, i.e., if $l-\mu>0$ is small. The case $l=\mu$ is  explicit; indeed, it is easy to check that the solution with initial value $v_0(x)=A\,|x|^{-\mu}$ is
\begin{equation}
v(x,t)=(Ce^{2(n-\mu)t}+1)^{1/(1-m)}|x|^{-\mu}, \quad C=A^{1-m}-1.
\end{equation}
For $A>1$ this solution blows up everywhere as $t\to\infty$ with rate $O(e^{2(n-\mu)t})$, while for $A<1$ it vanishes in finite time.

\medskip

\noindent {\bf 3.} The analysis of perturbations of the Barenblatt profiles,
$V_D$ with $D>0$, with large tails of the form  $v_0(x)-V_D(x)=O(|x|^{-l})$,
is an interesting related problem. The difference with the above analysis is
that the $v$-profile is regular, so no grow-up is expected if $l>\mu+2$.
Since the behaviour of $V_D$ at infinity is similar to the singular one,
$V_0$, and $V_D$ is still stationary, we also expect a continuum of
convergence rates  depending on $l$ from a certain range.  In this case we have
to mention that for $l>n$ there is a variational theory developed in the recent papers
\cite{BBDGV, BGV, PNAS}
that proves convergence with rate using the techniques of entropies, linearization
and functional inequalities.

\medskip

\noindent {\bf 4.} We could have used another of the possible scaling options, which is not adapted to the Barenblatt profiles but is still adapted to the singular solution. The simplest choice is
\begin{equation}
w(y,s)=[(1-m)(T-\tau)]^{-1/(1-m)}u(y,\tau), \quad s=(1-m) \log[ (T-\tau)/T]=(2/\beta)\, t,
\end{equation}
which leads to the equation
\begin{equation}
\frac{\partial w}{\partial s}= \Delta(w^m/m) +  w.
\end{equation}
Putting $w^m=Z$ and $p=1/m$ we get a variation of the Fujita equation
\begin{equation}
\frac{\partial Z^p}{\partial s}=a \Delta Z +b  Z^p.
\end{equation}
Studying this equation is equivalent to the study of the $v$ equation.
It is interesting to translate the results we have obtained and to
compare with the standard Fujita equation $u_t=\Delta u + u^p$.

\medskip

\noindent {\bf 5.} Our methods are not variational and our solutions do not belong to the usual spaces
of that theory, like spaces of finite relative energy or finite relative mass.

\vskip 1cm

\noindent {\large \sc Acknowledgment.}
Part of the work was done during visits of M.~F. and M.~W.
to the Univ. Aut\'onoma de Madrid and M.~F. to the Univ.
Duisburg--Essen. M.~F. was partially supported by the Slovak Research and Development
Agency under the contract No. APVV-0414-07 and by VEGA Grant 1/0465/09.
J.~L.~V. was partially supported by Spanish Project MTM2008-06326-C02.

%
%
%
%
%
%

\


  {\sc Addresses}

\noindent Marek Fila\\
 Department of Applied Mathematics and Statistics, Comenius University,\\
84248 Bratislava, Slovakia

\noindent Juan Luis V\'azquez\\
Departamento de Matem\'aticas, Universidad Aut\'onoma de Madrid,\\
28049 Madrid, Spain

and

\noindent Michael Winkler \\
Fakult\"at f\"ur Mathematik, Universit\"at Duisburg-Essen,\\
 45117 Essen, Germany\\
\

\vskip 1cm

  \noindent{\bf Key words:} Fast diffusion, extinction in finite time, extinction rate, nonlinear Fokker-Planck equation,  grow-up.

  \medskip

  \noindent {\bf MSC 2000:} 35K65, 35B40

\begin{thebibliography}{99}


\bibitem{BH} J. G. Berryman, C. J. Holland.
\textit{Stability of the separable solution for fast diffusion.}
Arch. Rat. Mech. Anal.
\textbf{74} (1980), 379--388.

\bibitem{BBDGV} {A. Blanchet, M. Bonforte, J. Dolbeault, G. Grillo, J. L. V\'azquez}.
{\em Asymptotics of the fast diffusion equation via entropy estimates,}
 Arch. Rat. Mech. Anal. {\bf 191} (2009), 347--385.

\bibitem{BGV}   M. Bonforte, J. Dolbeault, G. Grillo, J. L. V\'azquez.
\textit{Sharp rates of decay of solutions to the nonlinear fast diffusion equation via functional inequalities}, Preprint 2009, submitted.

\bibitem{PNAS}  {M. Bonforte, G. Grillo, J. L. V\'azquez}.
\textit{ Special fast diffusion with slow asymptotics. Entropy method
and flow on a Riemannian manifold},  Arch. Rat. Mech. Anal. (to appear),
published online: 1 July 2009; arXiv  0805.4750v1 [math.AP].


\bibitem{McD}
J. Denzler, R.~J. McCann.
\textit{Fast diffusion to self-Similarity: Complete spectrum, long-time asymptotics,
and numerology},  Arch. Rat. Mech. Anal. {\bf 175}  (2005), 301--342.


\bibitem{FS2000} E. Feiresl, F. Simondon.
\textit{Convergence for semilinear degenerate parabolic equations
in several space dimensions}, J. Dyn. Diff. Eq. \textbf{12} (2000),  647--673.

\bibitem{FKWY}
M. Fila, J. King, M. Winkler and E. Yanagida.
\textit{Optimal lower bound of the grow-up rate for a supercritical
parabolic equation}, J.~Diff. Equations {\bf 228} (2006), 339--356.

\bibitem{FW} M. Fila and M. Winkler.
\textit{Rate of convergence to a singular steady state
of a supercritical parabolic equation}, J.~Evol. Equations {\bf 8}
(2008), 673--692.

\bibitem{FWY1}
M. Fila, M. Winkler and E. Yanagida.
\textit{Grow-up rate of solutions for a supercritical
semilinear diffusion equation}, J.~Diff. Equations {\bf 205} (2004), 365-389.

\bibitem{Vsmooth}
J. L. V\'azquez.
{\em ``Smoothing and Decay Estimates for Nonlinear Diffusion
Equations''}, vol.~33 of Oxford Lecture Notes in Maths. and its
Applications, Oxford University Press, 2006.


\bibitem{Vbook} J. L. V\'azquez. {\em ``The
Porous Medium Equation. Mathematical Theory''}, Oxford
Mathematical Monographs, Oxford University Press, Oxford, 2007.


\end{thebibliography}
\end{document}